\documentclass[a4paper]{article}
\usepackage[comma,numbers,square,sort&compress]{natbib}
\renewcommand{\top}{{\mathrm{T}}}
\usepackage{graphicx,amsmath,amssymb}
\usepackage{subfigure,pgf,tikz,pst-node}
\usetikzlibrary{arrows,shapes,automata,positioning,fit}

\newtheorem{theorem}{Theorem}
\newtheorem{assumption}{Assumption}

\newtheorem{lemma}{Lemma}

\newtheorem{remark}{Remark}

\def\rm{\mathrm}
\begin{document}
\title{Coordination of Heterogeneous Nonlinear Multi-Agent Systems with Prescribed Behaviors}
\author{Yutao Tang \footnote{Y. Tang is with School of Automation, Beijing University of Posts and Telecommunications, Beijing 100876, China.  E-mail\,$: yttang@amss.ac.cn$}}
\date{}

\maketitle

{\noindent\bf Abstract} {In this paper, we consider a coordination problem for a class of heterogeneous nonlinear multi-agent systems with a prescribed input-output behavior which was represented by another input-driven system. In contrast to most existing multi-agent coordination results with an autonomous (virtual) leader, this formulation takes possible control inputs of the leader into consideration. First, the coordination was achieved by utilizing a group of distributed observers based on conventional assumptions of model matching problem. Then, a fully distributed adaptive extension was proposed without using the input of this input-output behavior. An example was given to verify their effectiveness.}


\section{Introduction}

In the past decades there has been a large percentage of multi-agent literature investigating on consensus-based coordination problem due to its numerous applications (see \cite{olfati2004consensus, ren2008distributed} and the references therein). Consensus means that a group of agents reach an agreement on a physical quantity of interest by interacting with their local neighbors.  Usually, a (virtual) leader is set up to define this quantity representing target trajectories or tasks. Including plenty of results for integrator-typed agents \cite{ren2008distributed, hong2006tracking, bauso2009consensus}, multi-agent systems with general linear dynamics have also be investigated even under variable topologies \cite{ni2010leader}. Recently, distributed/cooperative output regulation of multi-agent systems \cite{wang2010distributed, su2012cooperative, wang2014semi} was also proposed as a general framework for multi-agent coordination, which allows both reference tracking and disturbances rejection.

In most existing results, the quantity to be consensus on is assumed to be a constant or generated by an autonomous leader/exosystem, which may be restrictive or unpractical in some cases. On one hand, the leader might be a driven one especially when it is an uncooperative target or contains unmodeled uncertainties. Similar problems have been investigated by some authors in centralized or one-to-one cases \cite{saberi2001output}. On the other hand, the leader might be designed and tuned by us according to some objectives or from a high-level process, which is coincided with the classical model matching problem \cite{moore1972model, benedetto1986matching}. Also, a hierarchical control problem via abstraction was considered in \cite{girard2009hierarchical}, which was extended to a distributed version in \cite{tang2013hierarchical}, where the abstraction (virtual leader) is man-made with a tunable input in it. Therefore, it is necessary to study multi-agent control when the (virtual) leader contains driving inputs.

In fact, a coordinated tracking problem was investigated in \cite{hong2006tracking} when the agents' dynamics are integrators with inputs, an error bound was obtained by some input-to-state stability-like arguments. This problem was further analyzed in \cite{cao2012distributed} and the tracking error went to zero under a variable structure control law, which was extended to linear cases \cite{li2013tracking} by assuming that the leader and followers share the same dynamics. The case when the leader has different dynamics from that of those followers was later considered in \cite{tang2014leader}, where a disturbance decoupling condition was enforced to overcome the difficulties caused by heterogeneous dynamics. However, the verification of this geometry condition is nontrivial. Recently, a distributed generalized output regulation problem was formulated in \cite{tang2015auto} for a class of nonlinear agents and solved by an internal model-based controller combined with adaptive rules to deal with unknown-input leaders and possible (unbounded) disturbances.

Inspired by those works, we aim to investigate the coordination problem over a general class of heterogeneous nonlinear multi-agent systems with prescribed input-output behaviors described by an input-driven leader, whose dynamics is totally different from those of the other agents. The contribution of this work includes the following:
\begin{itemize}
	\item We extend the well-studied consensus problem \cite{hong2006tracking, ren2008distributed} to the case when the quantity to be consensus on is generated by an input-driven leader. While the leader here has a dynamics different from those of the non-identical followers, the results in \cite{cao2012distributed, li2013tracking} can be strictly recovered. When the leader has no driving inputs, these conclusions are consistent with the existing consensus or output regulation results \cite{ren2008distributed, su2012cooperative}.
	\item We extend the conventional model matching problem \cite{benedetto1986matching} and/or generalized output regulation formulation \cite{saberi2001output, ramos2004generalized} to their cooperative version for multi-agent systems with a driven leader.  In contrast to sufficient conditions and local results in \cite{ramos2004generalized} for the single agent case, we provide a necessary condition and global control laws for a class of nonlinear systems.  Additionally, the agents here are of much more general form and include many typical nonlinear systems, while only agents with unity relative degree was considered in \cite{tang2015auto}.
\end{itemize}

The rest of this paper is organized as follows. The problem is formulated in Section 2. Then our main results are presented in Section 3, followed by an example in Section 4. Finally, concluding remarks are presented in Section 5.

\section{Preliminaries and Problem Formulation}
Before the main results, we provide some preliminaries and then present the formulation of our problem.

\subsection{Graph theory and nonsmooth analysis}
Let $\mathbb{R}^n$ be the $n$-dimensional Euclidean space, $\mathbb{R}^{n\times m}$ be the set of $n\times m$ real matrices. $\text{diag}\{b_1,{\dots},b_n\}$ denotes an $n\times n$ diagonal matrix with diagonal elements $b_i\; (i=1,{\dots},n)$; $\text{col}(a_1,{\dots},a_n) = [a_1^{\top},{\dots},a_n^{\top}]^{\top}$ for any column vectors $a_i\; (i=1,{\dots},n)$.

A weighted directed graph (or weighted digraph) $\mathcal {G}=(\mathcal {N}, \mathcal {E}, \mathcal{A})$ is defined as follows, where $\mathcal{N}=\{1,{\dots},n\}$ is the set of nodes, $\mathcal {E}\subset \mathcal{N}\times \mathcal{N}$ is the set of edges, and $\mathcal{A}\in \mathbb{R}^{n\times n}$ is a weighted adjacency matrix . $(i,j)\in \mathcal{E}$ denotes an edge leaving from node $i$ and entering node $j$. The weighted adjacency matrix of this digraph $\mathcal {G}$ is described by $A=[a_{ij}]_{i,\,j=1,\dots,n}$, where $a_{ii}=0$ and $a_{ij}\geq 0$ ($a_{ij}>0$ if and only if there is an edge from agent $j$ to agent $i$). A path in graph $\mathcal {G}$ is an alternating sequence $i_{1}e_{1}i_{2}e_{2}{\cdots}e_{k-1}i_{k}$ of nodes $i_{l}$ and edges $e_{m}=(i_{m},i_{m+1})\in\mathcal {E}$ for $l=1,2,{\dots},k$. If there exists a path from node $i$ to node $j$ then node $i$ is said to be reachable from node $j$. The neighbor set of agent $i$ is defined as $\mathcal{N}_i=\{j\colon (j,i)\in \mathcal {E} \}$ for $i=1,...,n$. A graph is said to be undirected if $a_{ij}=a_{ji}$ ($i,j=1,{\dots},n$). The weighted Laplacian $L=[l_{ij}]\in \mathbb{R}^{n\times n}$ of graph $\mathcal{G}$ is defined as $l_{ii}=\sum_{j\neq i}a_{ij}$ and $l_{ij}=-a_{ij} (j\neq i)$. See \cite{mesbahi2010graph} for more details.

For the following nonsmooth analysis, we briefly review some basics of nonsmooth analysis. Consider the following differential equation with a discontinuous right-hand side:
\begin{align}\label{eq:non-smooth-sys}
\dot{x}=f(x,t)
\end{align}
where $f\colon \mathbb{R}^m\times \mathbb{R}\to \mathbb{R}^m$ is measurable and essentially locally bounded. A vector function $x(\cdot)$ is called a Filippov solution of \eqref{eq:non-smooth-sys} on $[t_0,\, t_1]$ if $x(\cdot)$ is absolutely continuous on $[t_0,\, t_1]$ and for almost all $t\in [t_0,\, t_1]$ satisfies the following differential inclusion: $\dot{x}\in \mathcal{F}[f](z,t)$, where $\mathcal{F}[f](z, t)=\bigcap_{\delta>0}\bigcap_{\mu(\bar N)=0}\bar{\mbox{co}}(f(B(z,\delta)-\bar N),t)$, $\bigcap_{\mu(\bar N)=0}$ denotes the intersection over all sets $\bar N$ of Lebesgue measure zero, $\bar{\mbox{co}}(E)$ is the convex closure of set $E$, and $B(z, \delta)$ denotes the open ball of radius $\delta$ centered at $z$.

Let $V\colon \mathbb{R}^m\to \mathbb{R}$ be a locally Lipschitz continuous function. The Clarke's generalized gradient of $V$ is defined by $\partial V(z)\triangleq\mbox{co}\{\lim_{i\to\infty}\nabla V(z_i): z_i\to z,\, z_i\notin \Omega_v\bigcup\bar N\}$, where $\mbox{co}$ denotes the convex hull, $\Omega_v$ is the set of Lebesgue measure zero where $\nabla V$ does not exist, and $\bar N$ is an arbitrary set of zero measure. The set-valued Lie derivative of $V$ with respect to \eqref{eq:non-smooth-sys} is defined as $\dot{\tilde V}\triangleq \bigcap_{\xi\in \partial V}\xi^T\mathcal{F}[f](z,t)$.

\subsection{Problem formulation}
Consider a group of heterogeneous nonlinear agents transformable into the following form
\begin{align}\label{sys:plant}
\begin{cases}
\dot{z}_i=A_i^0z_i+f_i(x_i),\\
\dot{x}_{i1}=x_{i2},\\
\vdots\\
\dot{x}_{in_i^x}=b_i^{\infty}u_i+g_{i}(z_i,x_i),\\
y_i=x_{i1}, \qquad i=1,{\dots},N
\end{cases}
\end{align}
where $z_i\in\mathbb{R}^{n_i^z}$, $x_i\triangleq\mbox{col}(x_{i1},\,\dots,\,x_{i n_i^x}) \in \mathbb{R}^{n_i^x}$, and $u_i,\, y_i\in \mathbb{R}$.  The matrix $A_i^0\in \mathbb{R}^{n_i^z\times n_i^z}$ is Hurwitz and the high-frequency gain $b_i^{\infty}$ is assumed to be a positive constant. With no loss of generality, we take $b_i^{\infty}=1$. The functions $f_i,\,g_i$ are infinitely differentiable.

The system \eqref{sys:plant} is general enough to cover some widely-investigated systems, including integrators, single-input single-output linear time-invariant systems and also some well-known nonlinear systems, e.g. controlled FitzHugh-Nagumo dynamics and controlled Van der Pol oscillator \cite{chen2015stabilization}.

We aim to drive all agents to match an input-output behavior described by
\begin{align}\label{sys:exosystem}
\dot{w}=Sw+dv(t),\quad y_{r}=c^{\top}w.
\end{align}
where $w\in \mathbb{R}^{n_0^w}$, $y_r\in \mathbb{R}$ and $v(t)\in \mathbb{R}$ is continuous satisfying $|v(t)|\leq l$ for some positive constant $l$. By the term ``matching", we mean to construct a proper (dynamic) controller such that the tracking error $e_i\triangleq y_i-y_{r}$ will asymptotically converge to zero for any $v$.  Since we only focus on input-output behavior of \eqref{sys:exosystem} which is often man-made according to some planning algorithms (e.g., by optimization), it is assumed, without loss of generality, to be minimal and has no zero dynamics. Thus, under a suitable coordinate transformation, we have
\begin{align*}
S&=\left[
\begin{array}{c|c}
0&I_{n_0^w-1}\\ \hline
s_0^0&s_1^0,\dots,s_{n-1}^0
\end{array}\right],\,\\
d&=\mbox{col}(0,\dots,0,d_{n_0^w}),\quad
c=\mbox{col}(1,0,\dots,0).
\end{align*}

In conventional model matching problem \cite{benedetto1986matching}, assuming the availability of $v$ and $w$ for all agents is reasonable and has been widely used, since the prescribed behavior is usually a mathematical description of the performance specifications. However, for multi-agent systems, we do not assume the availability of $w$ and $v$ to all agents to save resources, thus some agents may have no access to those information, which makes it much difficult to achieve collective behaviors.

To keep consistences, we denote the system \eqref{sys:exosystem} as a virtual leader in leader-following formulation \cite{hong2006tracking}. Associated with these multi-agent systems, a weighted digraph $\mathcal{G}$ can be defined with the nodes $\mathcal{N}=\{0,1,..., N\}$ to describe the communication topology, where node $0$ represent the leader. If the control $u_i$ can get access to the information of agent $j$, there is an edge $(j,\,i)$ in the graph $\mathcal{G}$, i.e., $a_{ij}>0$.  Also note that $a_{0i}=0$ for $i=1,...,N$, since the leader won't receive any information from the followers. Denote the induced subgraph associated with all followers as $\bar{\mathcal{G}}$.  A communication graph is said to be connected \cite{hong2006tracking} if the leader (node 0) is reachable from any other node of $\mathcal{G}$ and the induced subgraph $\bar{\mathcal{G}}$ of those followers is undirected. Given a communication graph $\mathcal{G}$, denote $H\in \mathbb{R}^{N\times N}$ as the submatrix of the Laplacian $L$ by deleting its first row and first column.

To achieve collective behaviors, the following assumption is often made.
\begin{assumption}\label{ass:graph}
	The communication graph $\mathcal{G}$ is connected.
\end{assumption}

Under this assumption, $H$ is positive definite by Lemma 3 in \cite{hong2006tracking}.  Denote its eigenvalues as $\lambda_1\geq\lambda_2\geq{\cdots}\geq \lambda_N>0$.

In this study, we mainly consider a control law $u_i$ of the form
\begin{align}\label{ctrl:dynamic}
\begin{split}
u_i&=\alpha_i(x_i,x^c_i,x_j,x^c_j,w)+\beta_i(x_i,x^c_i,x_j,x^c_j,w)v,\\
\dot{x}^c_i&=\zeta_i(x_i,x^c_i,x_j,x^c_j,w)+\gamma_i(x_i,x^c_i,w)v,\,j\in\mathcal{N}_i
\end{split}
\end{align}
where $x_i^c\in \mathbb{R}^{n_i^c}$ is a compensating variable for agent $i$ and the functions $\alpha_i(\cdot)$,\,$\beta_i(\cdot)$,\, $\zeta_i(\cdot)$,\, $\gamma_i(\cdot)$ are to be designed. When $n_i^c=0$, it reduces to a static control law. Here the input $v$ is some global information to all agents, and this control law becomes a distributed one when $v$ is not exactly used.

The coordination problem of multi-agent systems with a prescribed behavior can be formulated as follows. \emph{Given a multi-agent system composed of plant \eqref{sys:plant} and the behavior system \eqref{sys:exosystem}, find an integer $n_i^c$ and proper functions $\alpha_i(\cdot)$,\, $\beta_i(\cdot)$,\, $\zeta_i(\cdot)$,\, $\gamma_i(\cdot)$, such that, for any $v\in \mathbb{R}^{}$ and any initial condition $(z_i(0),\, x_i(0),\, x_i^c(0),\,w(0))\in \mathbb{R}^{n_i^z}\times \mathbb{R}^{n_i^x}\times \mathbb{R}^{n_i^c} \times \mathbb{R}^{n_0^w}$ of the composite system \eqref{sys:plant}-\eqref{ctrl:dynamic}, the tracking error $e_i\triangleq y_i-y_{r}$ satisfies
	\begin{align*}
	\lim\nolimits_{t\to+\infty} e_i=0, \,\mbox{ for any } i\in \{1,\,\dots,\,N\}.
	\end{align*}
}
\begin{remark}\label{rem:n=1}
	As having been mentioned, this problem resembles somehow the classical nonlinear model matching problem \cite{benedetto1986matching,isidori1995nonlinear} when $N=1$. Another related problem is so-called generalized output regulation \cite{saberi2001output, ramos2004generalized}, hence this work can be taken as their cooperative versions for nonlinear multi-agent systems, while only local results were presented in \cite{ramos2004generalized}, we provide here non-local controllers for nonlinear agents to solve it in a cooperative way.  When $v=0$ or $d_{n_0^w}=0$, this problem is exactly the existing consensus problem or the more general framework--distributed/cooperative output regulation--on a reference output which is assumed to be generated by an autonomous exosystem \cite{ren2008distributed, su2012cooperative}.
\end{remark}

The following assumptions on agents' dynamics are made to solve this problem.
\begin{assumption}\label{ass:relative-degree}
	For any $i\in \{1,\,\dots,\,N\}$, the relative degree of agent $i$ is no larger than that of the behavior system \eqref{sys:exosystem}, i.e., $n_i^x\leq n_i^w$.
\end{assumption}
\begin{assumption}\label{ass:lip} For any $i\in \{1,\,\dots,\,N\}$ and concerned $z_i,\,x_i$ , it holds $|g_i(z_i+\Delta,\,x_i)-g_i(z_i,\,x_i)|\leq M||\Delta||$ for some positive constant $M$.
\end{assumption}

\begin{remark}\label{rem:ass}
	The relative degree assumption \eqref{ass:relative-degree} is natural and can be proved to be necessary. In fact, it is also sufficient to solve our problem when $N=1$\cite[see][]{isidori1995nonlinear}.  Assumption \ref{ass:lip} is known as the Lipschitz property. When the concerned trajectories are contained in a compact set, this condition can be removed from the smoothness of $g_i(\cdot)$.
\end{remark}

\section{Cooperative Model Matching Design}
In this section, we constructively give control laws to solve the coordination problem of this multi-agent system with prescribed behaviors.  In conventional model matching literature \cite{benedetto1986matching, isidori1995nonlinear} and related generalized output regulation publications\cite{saberi2001output,ramos2004generalized}, the full information of system \eqref{sys:exosystem} is usually assumed to be available. However, this is not the case in multi-agent systems. While the availability of the control $v$ is possible, the information of $w$ might not be available for some agents.

To achieve collective behavior, the following distributed observer, inspired by the distributed design of \cite{hong2008distributed}, is employed for agent $i$ to estimate $w$ through the communication graph to facilitate our design.
\begin{align}\label{sys:observer-leader}
\dot{\eta}_i=S\eta_i+d v+ l_0c^{\top}\eta_{vi},
\end{align}
where $\eta_{vi}=\sum_{j=0}^N a_{ij}(\eta_i-\eta_j)$, $\eta_0=w$, $i=1,\,\dots,\,N$, and $l_0$ is a constant vector to be designed.  Letting $\bar \eta_i\triangleq\eta_i-w$ and denoting  $\bar \eta=\mbox{col}(\bar \eta_1,\dots,\bar \eta_N)$ gives
\begin{align}\label{sys:observer-leader-close}
\dot{\bar \eta}=[I_N \otimes S+H\otimes (l_0c^{\top})]\bar \eta.
\end{align}

The following lemma shows the effectiveness of these distributed observers.
\begin{lemma}\label{lem:observer-leader}
	Suppose $P$ is a positive definite matrix satisfying $S^\top P+P S<2Pdd^\top P$ and the communication graph is connected.  Taking $l_0=-\mu d^\top P$, there exists a constant $\mu^*$ such that when $\mu\geq \mu^*$, the system \eqref{sys:observer-leader-close} is uniformly exponentially stable in the sense of $||\bar \eta||\leq c_0  e^{-\lambda_0 t}$ for some positive constants $c_0$ and $\lambda_0$.
\end{lemma}
The proof is similar to Theorem 1 in \cite{ni2010leader} and omitted to save space.

\begin{remark}
	Note that a well-known sufficient condition for the solvability of the above linear matrix inequality is the stabilizability of $(S,\,d)$, thus this lemma holds naturally for our multi-agent systems. Although $\eta_0=w$ appears in \eqref{sys:observer-leader},  $y_{r}=c^{\top}\eta_0$ will suffice this design. The case when $y_{r}=w$ (i.e. the single-integrator case) has been partly considered in \cite{hong2006tracking,ren2008distributed}, while here the leader \eqref{sys:exosystem} contains an external input.
\end{remark}

With the help of those distributed observers for the followers, it is natural to replace $w$ by its estimation $\eta_i$. Now, we provide our first main theorem.

\begin{theorem}\label{thm:full-v}
	Under Assumptions \ref{ass:graph}--\ref{ass:lip}, the coordination problem of nonlinear multi-agent systems composed of \eqref{sys:plant} with a prescribed behavior described by system \eqref{sys:exosystem} is solvable by a control law of the form
	\begin{align}\label{ctr:fullstate+v}
	u_i&=-g_i(\xi_i,x_i)+x_{i(n_i^x+1)},\nonumber\\
	\dot{x}_{i(n_i^x+1)}&=x_{i(n_i^x+2)},\nonumber\\
	\vdots\nonumber\\
	\dot{x}_{in_0^w}&=\sum\nolimits_{j=1}^{n_0^w}s_j^0x_{ij}+d_{n_0}^wv+d_{n_0}^w\sum\nolimits_{j=1}^{n_0^w}k_j^0(x_{ij}-\eta_{ij}),\nonumber\\
	\dot{\xi}_i&=A_i^0\xi_i+f_i(x_i),\nonumber\\
	\dot{\eta}_i&=S\eta_i+d v+ l_0c^{\top}\eta_{vi}
	\end{align}
	where $l_0$ are defined in Lemma \ref{lem:observer-leader} and $k^0_1,\,\dots,\,k^0_{n_0^w}$ are selected constants such that the polynomial $s^{n_0^w}-k_{n_0^w}^0s^{n_0^w-1}-\dots-k_{1}^0$ is Hurwitz.
\end{theorem}
\noindent\textbf{Proof}. The proof will be split into two steps.

{\bf Step 1}: We first check the estimation performance of $\xi_i$ with respect to $z_i$. In fact, it can be found $z_i$-subsystem is in an output feedback form \cite{khalil2002nonlinear} with $x_i$ as its output. Letting $\bar z_i=\xi_i-z_i$ gives $\dot{\bar z}_i=A_i^0 \bar z_i$. Thus, $\xi_i$ will exponentially reproduce $z_i$ as times goes to infinity.

{\bf Step 2}: We now check the evolution of $\bar x_i=\hat x_i-w$, where $\hat x_i=\mbox{col}(x_{i1},\, \dots,\,x_{in_0^w})$. By \eqref{sys:plant}, \eqref{sys:exosystem}, \eqref{sys:observer-leader}, and \eqref{ctr:fullstate+v}, the composite system can be put into a compact form as follows.
\begin{align*}
\dot{z}_i&=A_i^0z_i+f_i(x_i),\\
\dot{\hat x}_{i}&=S\hat x_i+dv+dK^{\top}(\hat x_i-\eta_i)+\bar g_i(\bar z_i,\,z_i,\,x_i)\\
\dot{\eta}_i&=S\eta_i+d v+ l_0c^{\top}\eta_{vi}\\
\dot{\xi}_i&=A_i^0\xi_i+f_i(x_i)\\
\dot{w}&=Sw+dv
\end{align*}
Hence, one can obtain
\begin{align*}
\dot{\bar x}_{i}=(S+dK^{\top})\bar x_i-dK^{\top}\bar \eta_i+\bar g_i(\bar z_i,x_i)
\end{align*}
where $\bar g_i(\bar z_i,\,z_i,\, x_i)$ is a column vector function whose elements are zero except the $n_i^x$-th one. The $n_i^x$-th element is $g_i(z_i,\,x_i)-g_i(\xi_i,\,x_i)$.

Recalling that $\xi_i$ can exponentially reproduce $z_i$ and combining Assumption \ref{ass:lip} and Lemma \ref{lem:observer-leader}, we apply Lemmas 4.6 and 4.7 in \cite{khalil2002nonlinear} and obtain that $\bar x_i=0$ is an asymptotically stable equilibrium point of the $\bar x_i$-subsystem. Thus, $e_i$ will converge to 0 as $t\to+\infty$. Note that these arguments hold for any $i\in \{1,\,\dots,\,N\}$, thus we complete the proof. \hfill\rule{4pt}{8pt}

\begin{remark}
	Apparently, when $N=1$, this problem becomes to the conventional model matching problem \cite{benedetto1986matching} and/or generalized output regulation formulation \cite[refer to][]{saberi2001output, ramos2004generalized}. Hence, this formulation can be taken as their cooperative version for multi-agent systems with a driven leader.  In contrast to sufficient conditions and local results in \cite{ramos2004generalized} for the single agent case, we provide global control laws for a class of nonlinear systems.  Additionally, the agents here are of much more general form and include many typical nonlinear systems, while only agents with unity relative degree was considered in \cite{tang2015auto}.
\end{remark}

Since we merely have to match the input-output behavior of agent $0$, it might not be necessary to reproduce the full state $w$ of agent $0$. In fact, a reduced-order protocol will suffice our design as follows.

\begin{theorem}\label{thm:partial-v}
	Under Assumptions \ref{ass:graph}--\ref{ass:lip}, the coordination problem of nonlinear multi-agent systems composed of \eqref{sys:plant} with a prescribed behavior described by system \eqref{sys:exosystem} is solvable by a reduced-order control law of the form
	\begin{align}\label{ctr:partstate+v}
	u_i&=-g_i(\xi_i,x_i)+x_{i(n_i^x+1)},\nonumber\\
	\dot{x}_{i(n_i^x+1)}&=x_{i(n_i^x+2)},\nonumber\\
	\vdots\nonumber\\
	\dot{x}_{in_0^w}&=\sum\nolimits_{j=1}^{n_0^w}s_j^0x_{ij}+d_{n_0}^wv+d_{n_0}^wK\sum\nolimits_{j=0}^{N}a_{ij}(\hat x_{i}-\hat x_j),\nonumber\\
	\dot{\xi}_i&=A_i^0\xi_i+f_i(x_i)
	\end{align}
	where $\hat x_0=w$, $\hat x_i=\mbox{col}(x_{i1},\, \dots,\,x_{in_0^w})$,  $K=-\gamma d^T P$ with $\gamma\geq \frac{1}{\lambda_N}$, and $P$ is a positive definite matrix satisfying $S^\top P+P S<2Pdd^\top P$.
\end{theorem}
\noindent\textbf{Proof}. The proof is similar with that of Theorem \ref{thm:full-v}. After some mathematical manipulations, one can derive the composite system of agent $i$ in a compact form as follows.
\begin{align*}
\dot{z}_i&=A_i^0z_i+f_i(x_i),\\
\dot{\hat x}_{i}&=S\hat x_i+dv+dK\sum\nolimits_{j=0}^{N}a_{ij}(\hat x_{i}-\hat x_j)+\bar g_i(\bar z_i,\,z_i,\,x_i)\\
\dot{\xi}_i&=A_i^0\xi_i+f_i(x_i)\\
\dot{w}&=Sw+dv
\end{align*}
where $\bar g_i(\bar z_i,\,z_i,\, x_i)$ has been defined in the proof of Theorem \ref{thm:full-v}. Letting $\bar x_i=\hat x_i-w$ and $\bar z_i=\xi_i-z_i$ gives
\begin{align}\label{eq:error-thm2}
\begin{split}
\dot{\bar z}&=\bar A^0\bar z,\\
\dot{\bar x}&=(I_N\otimes S+ H\otimes dK)\bar x+\bar g(\bar z,\,z,\,x)
\end{split}
\end{align}
where $\bar A_0\triangleq \mbox{block}\,diag\{\bar A_1^0,\,\dots,\,\bar A_N^0\}$ and $\bar g(\bar z,\,z,\,x)$ is determined by $\bar g_i(\cdot)$.

Since $H$ is positive definite, there exists a unitary matrix $U$ such that $\Lambda \triangleq U^THU=\mbox{diag}\{\lambda_1,\dots,\lambda_N\}$. Let $J\triangleq(U^\top\otimes I_N)(I_N\otimes S+ H\otimes dK)(U\otimes I_N)$, then, $J=\mbox{block}diag\{J_1,\,\dots,\, J_N\}$ with $J_i\triangleq S+\lambda_i dK$. Apparently, $J_i^\top P+PJ_i=S^\top P+PS^\top-2\gamma\lambda_iPdd^\top P<0$. Hence $\bar A^1\triangleq I_N\otimes S+ H\otimes dK$ is Hurwitz. Recalling the stability of $\bar A^0$, there exist two positive definite matrices $\bar P^0$ and $\bar P^1$ satisfying $\bar {A^0}^\top\bar P^0+\bar P^0 \bar A^0=-I_{\sum_i n^z_i}$ and $\bar {A^1}^\top\bar P^1+\bar P^1 \bar A^1=-I_{\sum_i n^x_i}$.

We then take a quadratic Lyapunov candidate $\bar V=\varepsilon \bar z^\top\bar P^0\bar z+\bar x^\top \bar P^1\bar x$, where $\varepsilon>0$ is to be selected. It derivative along the trajectory of the above error system is then
\begin{align*}
\dot{\bar V}=-\varepsilon\bar z^\top \bar z-\bar x^\top \bar x+2\bar x^\top\bar P^1 \bar g(\bar z,\,z,\,x).
\end{align*}
This combined with the Lipschitzness of  $\bar g(\bar z,\,z,\,x)$ (Assumption \ref{ass:lip}) in $\bar z$ implies
\begin{align*}
\dot{\bar V}&\leq-\varepsilon\bar z^\top \bar z-\bar x^\top \bar x+2M||\bar P^1||||\bar x|| ||\bar z||\\
&\leq -(\varepsilon-2M^2||\bar P^1||^2)\bar z^\top \bar z-\frac{1}{2}\bar x^\top \bar x.
\end{align*}
Take $\varepsilon>2M^2||\bar P^1||^2+\frac{1}{2}$, then it follows for some positive constant $\bar \varepsilon$ that
\begin{align*}
\dot{\bar V}&\leq-\frac{1}{2}(\bar z^\top \bar z+\bar x^\top \bar x)\leq -\bar \varepsilon \bar V.
\end{align*}
which implies the convergence of $\bar x$, thus the conclusion is obtained. \hfill\rule{4pt}{8pt}

\begin{remark}
	It can be found that, when $v=0$ this problem reduces to the well-studied consensus problem, and these conclusions are consistent with the existing consensus or output regulation results \cite{hong2006tracking, ren2008distributed} when the quantity to be consensus on is generated by an autonomous leader. Additionally, while the leader here has a dynamics different from those of the non-identical followers, the results in \cite{cao2012distributed, li2013tracking} can be strictly recovered.
\end{remark}
\begin{remark}
	A similar problem has been considered in \cite{girard2009hierarchical} and \cite{tang2013hierarchical} under the formulation of hierarchical control, where the selected abstraction plays a similar role as prescribed behaviors in our formulation. The main difference between those two problems is that, we aim to achieve exactly tracking control, while a tradeoff was made in \cite{tang2013hierarchical} that we may sacrifice some accuracy without reallocating the designed controllers. Nevertheless, as the abstraction construction is still an open problem, these theorems might provide us a theoretical basis for abstraction selection to achieve better performances.
\end{remark}

\section{Fully Distributed Adaptive Extension}

In the last section, the cooperative controllers depend on the minimal eigenvalue $\lambda_N$ of $H$ and the leader's input $v$, which are actually global information. Usually, the multi-agent network is of a large scale and the eigenvalue is hard to compute. Also, the control input of \eqref{sys:exosystem} and even its upper bound might not be accessed by some agents. Thus, distributed control laws may be more favorable using only its local information.

Inspired by \cite{li2013tracking} and \cite{tang2014leader}, we propose an adaptive extension with non-smooth analysis to make proposed controllers fully distributed and achieve the coordination with prescribed behaviors. For simplicity, we only consider the reduced-order controller \eqref{ctr:partstate+v}.  With some modifications, we propose
\begin{align}\label{ctr:partstate-adaptive}
u_i&=-g_i(\xi_i,x_i)+x_{i(n_i^x+1)},\nonumber\\
\dot{x}_{i(n_i^x+1)}&=x_{i(n_i^x+2)},\nonumber\\
\vdots\nonumber\\
\dot{x}_{in_0^w}&=\sum\nolimits_{j=1}^{n_0^w}s_j^0x_{ij}-d_{n_0}^w\theta_id^\top P\hat x_{vi}-\theta_i\mbox{sgn}[(d^\top P\hat x_{vi})],\nonumber\\
\dot{\theta}_i&=||d^\top P\hat x_{vi}||_2^2+||d^\top P\hat x_{vi}||_1\nonumber\\
\dot{\xi}_i&=A_i^0\xi_i+f_i(x_i)
\end{align}
where $\hat x_{vi}=\sum\nolimits_{j=0}^{N}a_{ij}(\hat x_{i}-\hat x_j)$ and $\theta_i$ is the dynamic gain to be designed.

Note that the right-hand side of \eqref{ctr:partstate+v} is discontinuous, the stability of the closed-loop system will be analyzed by using differential inclusions and nonsmooth analysis \cite{cortes2008discontinuous}.  Since the sign function is measurable and locally bounded, by Proposition 3 in \cite{cortes2008discontinuous}, the Filippov solution of the closed-loop system exists. The following theorem shows the solvability of our problem by a fully distributed design.

\begin{theorem}\label{thm:partial-adaptive}
	Under Assumptions \ref{ass:graph}--\ref{ass:lip}, the coordination problem of nonlinear multi-agent systems composed of \eqref{sys:plant} with a prescribed behavior \eqref{sys:exosystem} is solvable by a distributed control of the form \eqref{ctr:partstate-adaptive}.
\end{theorem}
\noindent\textbf{Proof}. By some calculations, the composite system of agent $i$ can be put into a compact form.
\begin{align*}
\dot{z}_i&=A_i^0z_i+f_i(x_i),\\
\dot{\hat x}_{i}&=S\hat x_i-\theta_i d d^\top P\hat x_{vi}-\theta_i d\mbox{sgn}[(d^\top P\hat x_{vi})]+\bar g_i(\bar z_i,\,z_i,\,x_i)\\
\dot{\theta}_i&=||d^\top P\hat x_{vi}||_2^2+||d^\top P\hat x_{vi}||_1\\
\dot{\xi}_i&=A_i^0\xi_i+f_i(x_i)\\
\dot{w}&=Sw+dv
\end{align*}
Performing a coordinate transformation $\bar x_i=\hat x_i-w$ and $\bar z_i=\xi_i-z_i$ gives
\begin{align}\label{eq:thm3-closed}
\dot{\bar z}_i&=A_i^0\bar z_i,\nonumber\\
\dot{\bar x}_{i}&=S\bar x_i-\theta_idd^\top P\bar x_{vi}-d[\theta_i \mbox{sgn}[(d^\top P\bar x_{vi})]+v]+\bar g_i(\bar z_i,\,z_i,\,x_i)\nonumber\\
\dot{\theta}_i&=||d^\top P\hat x_{vi}||_2^2+||d^\top P\hat x_{vi}||_1
\end{align}
where $\bar x_{vi}=\sum\nolimits_{j=0}^{N}a_{ij}(\bar x_{i}-\bar x_j)$ and $\bar x_0=0$. Also note that
\begin{align*}
\dot{\bar z}&=\bar A^0\bar z,\\
\dot{\bar x}&=(I_N\otimes S-EH\otimes dd^\top P)\bar x-(E\otimes d)\mbox{sgn}[(H\otimes d^\top P)\bar x]-(I_N\otimes d) v+\bar g(\bar z,\,z,\,x)
\end{align*}
where $\bar A_0$, $\bar g(\bar z,\,z,\,x)$ are defined as above in Equation \eqref{eq:error-thm2} of Theorem \ref{thm:partial-v}, $E=\mbox{diag}(\theta_1,\,\dots,\,\theta_N)$, and $\mbox{sgn}[(H\otimes d^\top P)\bar x]$ is defined elementwise.

Since the matrix $\bar A_0$ is Hurwitz by assumptions, there exists a positive definite matrix $Q$ such that $\bar A_0^\top Q+Q\bar A_0^\top=-I_{n_z}$, where $n^z=\sum_i n_{i}^z$.  To prove this theorem, we consider a Lyapunov candidate as follows.
\begin{align}
V=\bar x^\top (H\otimes P)\bar x+ \kappa \bar z^\top Q\bar z+\sum\nolimits_{i=1}^N(\theta_i-\Theta)^2
\end{align}
where $\kappa$ and $\Theta$ are positive constants to be designed. Its set-valued Lie derivative along the trajectory of the closed-loop system under controller \eqref{ctr:partstate-adaptive} is
\begin{align*}
\dot{\tilde V}&\triangleq\mathcal{F}[\bar x^\top (H\otimes P)\dot{\bar x}-\kappa \bar z^\top \bar z+2\sum\nolimits_{i=1}^N(\theta_i-\Theta)\dot{\theta}_i]\\
&=2\bar x^\top (H\otimes P)(I_N\otimes S-EH\otimes dd^\top P)\bar x-2\bar x^\top (H\otimes P)(E\otimes d)\mbox{sgn}[(H\otimes d^\top P)\bar x]\\
&-2\bar x^\top (H\otimes P)(I_N\otimes d) v+2\bar x^\top (H\otimes P)\bar g(\bar z,\,z,\,x)]-\kappa \bar z^\top \bar z+2\sum\nolimits_{i=1}^N(\theta_i-\Theta)\dot{\theta}_i\\
&=\Gamma_1+\Gamma_2+2\bar M||\bar x||_2||z||_2-\kappa \bar z^\top \bar z+2\sum\nolimits_{i=1}^N(\theta_i-\Theta)\dot{\theta}_i
\end{align*}
where $\Gamma_1\triangleq\bar x^\top [H\otimes (PS+S^\top P) ]\bar x-2\bar x^\top (HEH\otimes Pdd^\top P)\bar x$,\, $\Gamma_2\triangleq-2\bar x^\top (HE\otimes Pd)\mbox{sgn}[(H\otimes d^\top P)\bar x]-2\bar x^\top (H\otimes Pd) v$, and $\bar M$ is a constant greater than $M||H\otimes P||_2$. Here we use the fact that $F[f](x)={f(x)}$ if $f(x)$ is continuous \cite{cortes2008discontinuous} since $2\bar x^\top (HE\otimes Pd)\mbox{sgn}[(H\otimes d^\top P)\bar x] =\sum_{i=1}^N \theta_i||d^\top P \bar x_{vi}||_1$ and $\Gamma_2$ is continuous.

Let $\xi=(U^\top \otimes I_N)\bar x$ and $\hat x_{vi}=\sum\nolimits_{j=0}^{N}a_{ij}(\hat x_{i}-\hat x_j)$, where $U$ is defined in Theorem \ref{thm:partial-v}, then $\hat x_v=(H\otimes I_N)\bar x$. One  can obtain
\begin{align*}
\Gamma_1&=\xi^\top[\Lambda\otimes (PS+S^\top P) ]\xi-2\Theta \xi^\top (\Lambda^2 \otimes Pdd^\top P)\xi- 2\bar x^\top [H(E-\Theta I_N)H\otimes Pdd^\top P ]\bar x\\
&=\sum_i\lambda_i\xi_i^\top[PS+S^\top P-2\Theta\lambda_i Pdd^\top P)\xi_i-2\sum_i (\theta_i-\Theta)\hat x_{vi}^\top Pdd^\top P\hat x_{vi}
\end{align*}
Let $\Theta$ be large enough such that $\Theta\lambda_i>1$, then there exists a positive constant $\varepsilon$ such that
\begin{align}\label{eq:thm3-1}
\Gamma_1\leq  -\varepsilon\bar x^\top \bar x-2\sum_i (\theta_i-\Theta)\hat x_{vi}^\top Pdd^\top P\hat x_{vi}
\end{align}
For the second term, we have
\begin{align}\label{eq:thm3-2}
\Gamma_2&=-2\sum_i\Theta ||d^\top P\hat x_{vi}||_1-2\bar x^\top (H\otimes Pd) v-2\bar x^\top [H(E-\Theta I_N)\otimes Pd]\mbox{sgn}[(H\otimes d^\top P)\bar x]\nonumber\\
&\leq -2\sum_i(\Theta-l) ||d^\top P\hat x_{vi}||_1-2\sum_i(\theta_i-\Theta) ||d^\top P\hat x_{vi}||_1
\end{align}
Note that by Young's inequality we have
\begin{align}\label{eq:thm3-3}
2\bar M||\bar x||_2||z||_2\leq \frac{\varepsilon}{2} ||\bar x||^2+ \frac{2\bar M^2}{\varepsilon}||z||^2
\end{align}
By letting $\Theta=\max \{l,\, \frac{1}{\lambda_N}\}$ and $\kappa=\frac{2\bar M^2}{\varepsilon}+1$ and combining \eqref{eq:thm3-1}-\eqref{eq:thm3-3}, we have
\begin{align}\label{eq:thm3-lyapunov}
\dot{V}&\leq -\frac{\varepsilon}{2}\bar x^\top \bar x-2\sum_i(\Theta-l) ||d^\top P\hat x_{vi}||_1-\bar z^\top \bar z+2\sum\nolimits_{i=1}^N(\theta_i-\Theta)(\dot{\theta}_i-\Delta_i)\nonumber\\
&\leq -\frac{\varepsilon}{2}\bar x^\top \bar x-\bar z^\top \bar z\triangleq W(\bar x, \,\bar z)
\end{align}
where  $\Delta_i=||d^\top P\hat x_{vi}||_2^2+||d^\top P\hat x_{vi}||_1$.

Apparently, the trajectory of the closed-loop system is bounded and thus the derivatives of $\bar x$ and $\bar z$ are also bounded from \eqref{eq:thm3-closed}. Hence, $W$ is uniformly continuous with respect to the time $t$. By integrating the both sides of \eqref{eq:thm3-lyapunov}, we have
$$
\int_{t_0}^{\infty}W(\bar x(t),\,\bar z(t))dt\leq V(t_0).
$$
Recalling Barbalat's lemma \cite{khalil2002nonlinear}, we have $W(\bar x,\,\bar z)\to 0$ when $t\to \infty$, and hence $\bar x$ converges to zero when $t$ goes to infinity, while $\theta_i$ converges to some finite value. \hfill\rule{4pt}{8pt}
\begin{remark}\label{rem:li-tang}
	Although similar design has been used in \cite{cao2012distributed, li2013tracking}, the agents considered here are of much general heterogeneous dynamics and include existing results as its special cases. Furthermore, it can also be proved that the disturbance decoupling condition used in \cite{tang2014leader} is sufficient for Assumption \ref{ass:relative-degree}, while the latter one is checkable.
\end{remark}
\begin{remark}
	Since the discontinuous signum function is employed in this fully distributed design, the unfavorable chattering might rise and result in some instability of this control law. We thus propose one of its continuous approximation as follows:
	\begin{align}\label{ctr:partstate-adaptive-sigma}
	u&=-g_i(\xi_i,x_i)+x_{i(n_i^x+1)},\nonumber\\
	\dot{x}_{i(n_i^x+1)}&=x_{i(n_i^x+2)},\nonumber\\
	\vdots\nonumber\\
	\dot{x}_{in_0^w}&=\sum\nolimits_{j=1}^{n_0^w}s_j^0x_{ij}-d_{n_0}^w\theta_id^\top P\hat x_{vi}-\theta_i\mbox{sat}_\epsilon(d^\top P\hat x_{vi}),\nonumber\\
	\dot{\theta}_i&=||d^\top P\hat x_{vi}||_2^2+||d^\top P\hat x_{vi}||_1-\sigma\theta_i,\nonumber\\
	\dot{\xi}_i&=A_i^0\xi_i+f_i(x_i)
	\end{align}
	where ${\rm sat}_\epsilon(x)=
	\begin{cases}
	x/\epsilon, &\text{ if } |x|\leq\epsilon;\\
	{\rm sgn}(x/\epsilon),& \text{ if } |x|>\epsilon
	\end{cases}$, and $\sigma, \epsilon>0$ are tunable parameters.  It can be verified following a similar proof as that in \cite{jiang1999robust} and \cite{tang2015auto} that this control law will eventually drive all tracking errors into a bounded set. Furthermore, the bound can be sufficiently small by tuning $\sigma$ and $\epsilon$ according to practical control goals.
\end{remark}

\section{Examples}

We give an example to illustrate the effectiveness of our design in previous sections.

Consider three nonlinear agents including a controlled damping oscillator, a controlled FitzHugh-Nagumo dynamics and a controlled Van der Pol oscillator \cite{chen2015stabilization} as follows.
\begin{align*}
\mbox{Agent}~1:\,
\begin{cases}
\dot{x}_{11}=x_{12}\\
\dot{x}_{12}=\Delta_1+u_1\\
y_1=x_{11}
\end{cases}\quad
\mbox{Agent}~2:\,
\begin{cases}
\dot{z}_{21}=-c_2z_{21}+b_2x_{21}\\
\dot{x}_{21}=\Delta_2+u_2\\
y_{2}=x_{21}, \, c_2>0
\end{cases}\quad
\mbox{Agent}~3:\,
\begin{cases}
\dot{x}_{31}=x_{32}\\
\dot{x}_{32}=\Delta_3+u_3\\
y_3=x_{31}.
\end{cases}
\end{align*}
where $\Delta_1=-x_{11}-x_{12}$, $\Delta_2=x_{21}(a_2-x_{21})(x_{21}-1)-z_{21}$, and $\Delta_3=-x_{31}+a_3 (1-x_{31}^2)x_{32}$.  The prescribed input-output behavior (agent 0) is $\dot{w}_1=w_2,\;\dot{w}_2=v,\;y_r=w_1$.

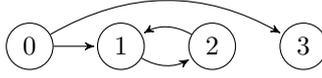
\begin{figure}
	\centering
	\begin{tikzpicture}[shorten >=1pt, node distance=1.2 cm, >=stealth',
	every state/.style ={circle, minimum width=0.2cm, minimum height=0.2cm}, auto]
	\node[align=center,state](node0) {0};
	\node[align=center,state](node1)[right of=node0]{1};
	\node[align=center,state](node2)[right of=node1]{2};
	\node[align=center,state](node3)[right of=node2]{3};
	\path[->]   (node0) edge (node1)
	(node1) edge [bend right] (node2)
	(node2) edge [bend right] (node1)
	(node0) edge [bend left]  (node3)
	;
	\end{tikzpicture}
	
	\caption{The communication graph .}
	\label{fig:graph}
\end{figure}

\begin{figure}
	\centering
	\subfigure[]{
		\includegraphics[width=0.4\textwidth]{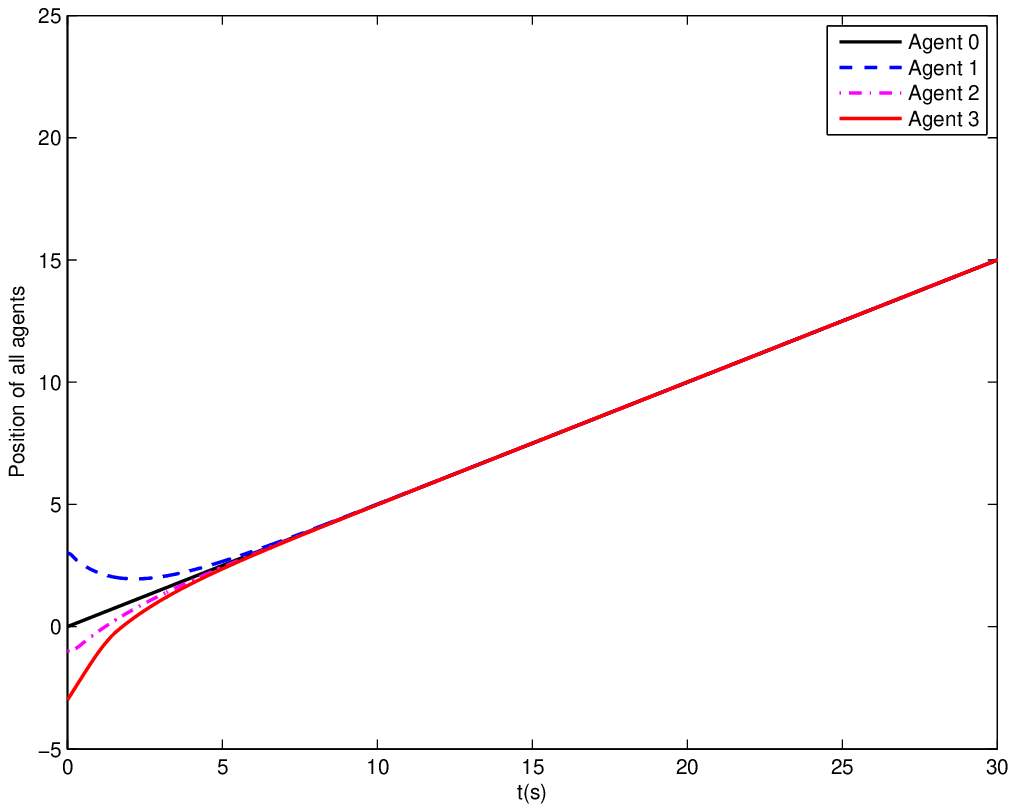}
	}
	\subfigure[]{
		\includegraphics[width=0.4\textwidth]{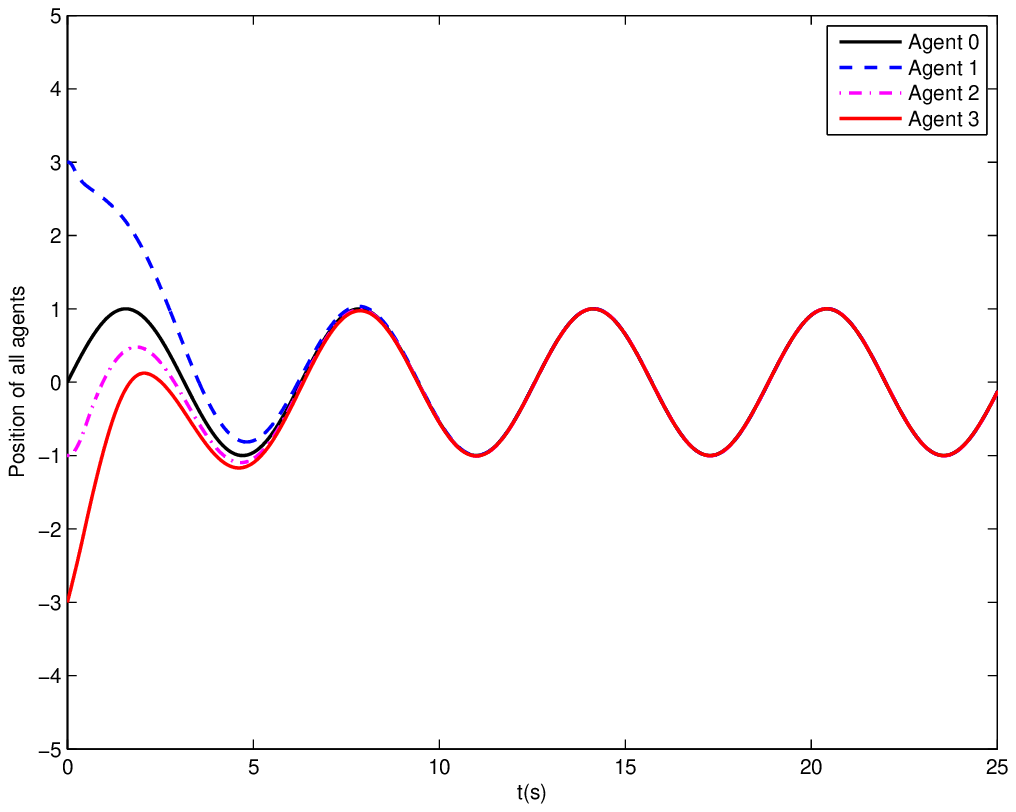}
	}
	\caption{\footnotesize The performance of the controller \eqref{ctr:partstate-adaptive}: Figure a) for $v=0$ and Fig. b) for $v=-w_1$.}\label{fig:simu}
\end{figure}

It can be verified that the communication topology in Fig.~\ref{fig:graph} is connected and Assumptions 1 and 2 are also satisfied. We employ the distributed control \eqref{ctr:partstate-adaptive} to solve this problem.  Take $a_2=1,\,b_2=1,\,c_2=1,\, a_3=1$ and the initial value of each state variable randomly generated during $[-3,\,3]$. First, we let $v=0$ to generate a ramping signal, and then take $v=-w_1$ to generate a sinusoidal signal. While the controller is fixed, the simulation results are depicted in Fig.~\ref{fig:simu}.

\section{Conclusions}
A coordination problem with prescribed behaviors for a class of nonlinear heterogeneous multi-agent systems was formulated as a distributed extension of conventional model matching problem. Based on some conventional assumptions, two control laws and a fully distributed extension were given to solve this problem. Future works include nonlinear MIMO multi-agent systems and with more general graphs.

\bibliographystyle{IEEEtr}
\bibliography{mybib}

\end{document}